\documentclass[12pt]{article}

\title{Diophantine Approximation on varieties II: Explicit estimates for Arithmetic Hilbert Functions} %and arithmetic Interpolation}%, arithmetic Regularity and Interpolation}

\author{Heinrich Massold}

\usepackage{bezier,amsfonts}

\usepackage{amsmath,amsfonts,bbm,amssymb}

\newtheorem{Satz}{}[section]

\newcommand{\satz}[1]{\vspace{2mm} \begin{Satz}{\bf #1}}

\newcommand{\proof}{\vspace{3mm} {\sc Proof}\hspace{0.2cm}}

\newcommand{\la}{\langle}

\newcommand{\ra}{\rangle}

\newcommand{\di}{{\mbox{div}}}

\newcommand{\rk}{\mbox{rk}}

\newcommand{\R}{{\mathbbm R}}

\newcommand{\Pe}{{\mathbbm P}}

\newcommand{\Z}{\mathbbm{Z}}

\newcommand{\N}{\mathbbm{N}}

\newcommand{\CH}{{\cal H}}

\newcommand{\CB}{{\cal B}}

\newcommand{\CO}{{\cal O}}

\newcommand{\CX}{{\cal X}}

\newcommand{\CY}{{\cal Y}}

\newcommand{\CZ}{{\cal Z}}

\newcommand{\CL}{{\cal L}}

\newcommand{\CU}{{\cal U}}

\newcommand{\fa}{{\mathfrak a}}

\newcommand{\fP}{{\mathfrak P}}

\newcommand{\Q}{{\mathbbm Q}}

\newcommand{\C}{{\mathbbm C}}

\newcommand{\A}{{\mathbbm A}}

\newcommand{\wde}{\widehat{\deg}}

\newcommand{\spec}{\mbox{Spec}}

\begin{document}

\parindent0mm

\maketitle

\thispagestyle{empty}

\tableofcontents

\section{Introduction}

Let $X$ be a projective variety over a field $k$, i.e.\@ a reduced projective 
scheme over $\spec \; k$, and $\CX$
a projective variety of relative dimension $s$ over $\Z$ 
that is a reduced projective flat scheme over $\spec \Z$. For an ample line 
bundle $L$ on $X$, a positive metrized line bundle $\bar{\CL}$ 
on $\CX$ and a $D \in \N$ 
defines $H_X(D)$ as the dimension of the vector space of global sections
of $L^{\otimes D}$ on $X$, and 
$\hat{\CH}_{\CX}(D)$ as the arithmetic degree of the arithmetic
bundle of global sections of $\bar{\CL}^{\otimes D}$ on $\CX$. In this context,
there are the well know algebraic and arithmetic Hilbert-Samuel formulas
\[ H_X(D) = \deg_L X \frac{D^s}{s!} + O(D^{s-1}), \]
\[ \hat{\CH}_\CX(D) = h_{\bar{\CL}} (\CX) \frac{D^{s+1}}{(s+1)!}
   + O(D^s \log D). \]
These formulas allow to calculate the Hilbert functions as $D$ goes
to infinity, but do not give any explicit estimates for a given $D$. Therefore,
they are rather useless in the context of diophantine approximations.
For subvarieties of projective space $\Pe^t$ the explicit upper bounds 
\[ H_X(D) \leq \deg Y {D +s \choose s}, \]
\[ \hat{H}_\CX(D) \leq \left(Dh(\CY) +c_1 \deg X D + \deg X
    \left(\frac12 \log \deg X+s\log D\right)\right) {D+s \choose s} \]
were proved in \cite{Ch},\cite{LNM1752}, chapter 9.1 and \cite{CP}.
Corresponding lower bounds for the Hilbert functions do not hold in 
general, but only for $D$ bigger than a number $\bar{D}$
depending on $\CX$ and $X$ respectively. 
For irreducible $\CX,X$ consider the weak estimates
\[ H_X(D) \geq c \deg X D^s, \quad \mbox{and} \quad
   \hat{H}_{\CX}(D) \geq (c'h(\CX) + c'' \deg X) D^{s+1}, \]
where $c,c',c''$ are supposed to be positive constants only depending on
$s$ and $t$, and for general $\CX,X$ the strong estimates
\[ H_X(D) \geq  \deg X {D-(t-s)\bar{D}+s \choose s}, \]
and
\[ \hat{H}_{\CX}(D) \geq (h(\CX)- c'' \deg X) 
                D {D-(t-s)\bar{D}+s \choose s}. \]
This paper contains proofs for the weak estimates. The proofs use a 
relatively simple induction argument. The stronger estimates require
the concepts of regularity and conjecturaly arithmetic regularity.

The concepts and results of the first part of this series on
diophantine approximation will
not be used in this paper though there are similar prerequisites
and argumentations.

\section{Locally complete intersections and algebraic \\ Hilbert functions}

In this section $k$ denotes an arbitrary field, and 
$\Pe^t=\mbox{Proj}(k[x_0,\ldots,x_t])$ projective
space over $k$. A closed subset $X \subset \Pe^t$ is referred to as
projective subvariety.

%Let $A$ be a graded $\CO_k$ algebebra, and $\CM$ a finitely generated
%graded $A$ module. For any tupel ${\bf a} = (a_1, \ldots, a_s) \in A^s$ let
%$H_{s-i}({\bf a},\CM)$ be the homology groups of the Koszul comples
%$K_\bullet ({\bf a},\CM)$. They inherit a graduation from $M$.
%Further denote $A_k = A \otimes_{\CO_k} k$, and
%$M = \CM_k = M_{\otimes \CO_k} k$.

%\satz{Definition} \label{dr1}
%Let $D \in \N$, and $\CI \subset A$ be a homogenous ideal with base
%extension $I = \CI_k$. 
%$M$ is called $(D,I)$-regular if for every
%tupel ${\bf a} = (a_1, \ldots, a_s) \in A_k^s$ such that $\dim M/IM = 0$,
%\[ (I H_{s-i}({\bf a},M))_D = 0, \quad \forall i, \quad \forall \nu > D-i. \]
%A subvariety $X \subset \Pe^t_k$ is called locally $D$-regualar if
%for some nonzero 
%homogenous ideal $I \subset k[x_0,\ldots,x_t]$ the homogeneous
%oordinate ring $A_k$ of $X$ is $(D,I)$-regular.
%\end{Satz}

%\satz{Lemma}
%For a zero dimensional $D$- regular variety $X$, the vanishing
%ideal $I_X$ is generated by elements of degree at most $D$.
%\end{Satz}

\satz{Definition}
A subvariety $X$ of $\Pe^t_K$ of pure codimension $s$ is called a locally 
complete intersection of hypersurfaces $H_1, \ldots H_s$ if there is an
open subset $U \subset \Pe^t$ such that
\begin{equation} \label{lvd}
X = \overline{H_1 \cap \cdots \cap H_s \cap U}.
\end{equation}

Let $Y$ be an irreducible subvariety of $\Pe^t$.
A subvariety $X$ of pure codimension $s$ is called a locally complete
intersection at $Y$ if there are hypersurfaces $H_1, \ldots, H_s$ such that
$X$ is the union of the irreducible
components of $H_1 \cap \cdots \cap H_s$ that contain $Y$.

\end{Satz}

\satz{Lemma} \label{bas}

\begin{enumerate}

\item
If $X \subset \Pe^t$ is a locally complete intersection, and 
$Y \subset X$ is a subvariety of codimension zero, then $Y$ is a
locally complete intersection.

\item
If $Y$ is an irreducible variety and $X$ is a locally complete intersection
at $Y$, then $X$ is a locally complete intersection.

\item
If $X$ is a locally complete intersection at $Y$, and $Z$ a subvariety
that contains $Y$ and intersects $X$ properly, the union $W$
of the components of $X \cap Z$  that contain $Y$ is a locally complete 
intersection at $Y$.
\end{enumerate}
\end{Satz}

\proof
1. In a representation of $X$ as locally complete intersection, one 
only needs to replace $U$ by 
$U' = U \setminus \overline{X\setminus Y}$.

\vspace{2mm}

2. If $X$ is the union of irreducible components of $H_1 \cap \cdots \cap H_s$
that contain $Y$, let $\bar{X}$ be the union of irreducible components
of $H_1 \cap \cdots \cap H_s$ that do not contain $Y$, and
$\CU = \Pe^t \setminus \bar{X}$. Then,
$X= \overline{U \cap H_1 \cap \cdots \cap H_s}$.

\vspace{2mm}

3. Is obvious.

\satz{Proposition} \label{alHilbert}
Let $X \subset \Pe^t$ be a subvariety of pure dimension $s$ in $\Pe^t$, 
and denote by 
\[ H_X(D) = \dim H^0(X,O(D)) \]
the algebraic Hilbert function.

\begin{enumerate}
\item
For every $D \in \N$,
\[ H_X(D) \leq \deg X {D+t-s \choose t-s}. \]

\item
If $X$ is a locally complete intersection of hypersurfaces of degree
$D_1, \ldots,D_s$, then for $D\geq \bar{D}:=D_1+\cdots+D_s-s$,
\[ H_X(D) \geq \deg X {D- \bar{D} + t-s  \choose t-s}. \]
\end{enumerate}
\end{Satz}

\proof
1. \cite{Ch}, Theorem 1.

\vspace{1mm}

2. \cite{CP}, Corollaire 3.

\satz{Lemma} \label{summe}
Let $k$ be of characteristic zero,
$n \in \N$, and $v_1, \ldots, v_n \in k^n$ such that for each
$i=1,\ldots,n$ the $i$th component of $v_i$ is nonzero. Then, there are
$m_1,\ldots,m_n \in \N$ with $m_i \leq n$ such that in 
\[ v = m_1 v_1 + \cdots + m_n v_n \in k^n \]
no component is zero.
\end{Satz}

\proof
The Lemma clearly holds for $n=1$; so assume it holds for $n-1$.
Then, for $v_1, \ldots, v_n$ as in the Lemma, there are
$m_1, \ldots, m_{n-1}$ with $m_i \leq n-1, i=1,\ldots,n-1$ such that
\[ w = m_1 v_1 + \cdots + m_{n-1} v_{n-1} \]
has the first $n-1$ components not equal to zero. 
If also the last component of $w$ is nonzero, choosing $m_n=0$ proves the
Lemma. If the $n$th component of $w$ equals zero, let $u = v_n$, and 
define $w_j$ as the $j$th component of $w$, and $u_j$ as the $j$th
component of $u$. Further, $k_j := u_j/w_j, j=1, \ldots, n-1$.
As these are $n-1$ numbers there is an $m_n \neq 0$ with $m_n \leq n$ such
that $\frac1{m_n} \neq - k_j$ for all $j=1,\ldots,n-1$. 
Consequently, for $v_j$ the
$j$the component of $v = w + m_nu$, we have $v_j = w_j + m_n u_j$.
Thus, for $j\leq n-1$, if $u_j=0$, then $v_j = w_j \neq 0$, and if
$v_j \neq w_j$, then $m_n \neq -k_j \neq 0$, and consequently
$v_j = w_j + m_n u_j \neq w_j -\frac1{k_j} u_j =0$.
Since $w_n=0$, also $v_n = m_n u_n \neq 0$, and
thus $v_j \neq 0 $ for all $j=1,\ldots,n$.

\satz{Corollary} \label{summe1}
Let $X,Y \subset \Pe^t$ be algebraic varieties over 
a field $k$ of characteristic zero,
and $X= X_1 \cup \cdots \cup X_n$ the decomposition into irreducible
components. If for some $D \in \N$ and every $i=1,\ldots,n$, there is an
$f_i \in H^0(Y,O(D))$ that is zero on $Y$ and nonzero on $X_i$, then there
is an $f \in H^0(Y,O(D))$ such that $f$ is zero on $Y$ and nonzero on 
every $X_i$.
\end{Satz}

\proof
For $i=1,\ldots,n$,
let $x_i \in X_i(\bar{k})$ such that $(f_i)_{x_i} \neq 0$, and
$v_i \in \bar{k}$ the vector whose $j$th component is $(f_i)_{x_j}$. 
By the Lemma, there are $m_1,\ldots,m_n \leq n$ such that  
with $f=m_1f_1 + \cdots +m_nf_n$  one has
$f_{x_i} \neq 0$ for every $i$, hence the restriction of $f$ to every $X_i$ is
nonzero.

\satz{Corollary and Definition} \label{TR}
In characteristic zero,
every pure dimensional subvariety $X$ of $\Pe^t$ is a subvariety of
codimension $0$ of a complete intersection, and consequently a locally complete
intersection. Define the top regularity
$\bar{D}=\mbox{TR}(X)$ of $X$ 
as the minimal number $\bar{D} \in \N$ such that $X$ is a locally
complete intersection of hypersurfaces of degree at most $\bar{D}$.

If $s$ is the codimension of $X$, and $\bar{X}$ a subvariety of codimension 
$r \leq s$ containing $X$, then $\deg X \leq \deg \bar{X} \;(TR(X))^{s-r}$.
Further, if $X,Y,Z,W$ are as in Lemma \ref{bas}.3, then
$TR(W) \leq \max(TR(X),TR(Z))$.
\end{Satz}

\proof
Let $s$ be the codimension of $X$, and assume the claim holds for subvarieties
of codimension at most $s-1$. Further, let $X_1$ be a subvariety of
$\Pe^t$ of codimension $s-1$ that contains $X$. By induction hypothesis,
there are homogeneous polynomials $h_1,\ldots,h_{s-1}$ such that
$X_2=V(h_1) \cap \cdots \cap V(h_1)$ has codimension $s-1$ and contains
$X_1$. By the previous corollary, there is a homogeneous polynomial $h_s$ 
that is zero on $X$ and nonzero on every irreducible component of
$X_2$. Hence $V(h_1) \cap \cdots \cap V(h_s)$ has codimension $s$ and
contains $X$. To see that $X$ is a locally complete intersection, 
let $X_3$ be the union of irreducible components of
$V(h_1) \cap \cdots \cap V(h_s)$ except for $X$, and choose
$U = \Pe^t \setminus X_3$. Then,
$X= \overline{U \cap V(h_1) \cap \cdots \cap V(h_s)}$, proving the first claim.

Since there are homogeneous polynomials 
$h_1$ with $\deg h_i \leq \bar{D}$, and an open set $U$ with
$X = \overline{U \cap V(h_1) \cap \ldots \cap V(h_s)}$, and
$s>r$, for every irreducible component $Z$ of $\bar{X}$ there is an
$h_i$ with $(h_i)|_Z \neq 0$. Hence, by Lemma \ref{summe}, there is
a linear combination $h$ of the $h_i$ with coefficients in $k$ that
has nonzero restriction to every irreducible component of $\bar{X}$,
and consequently $\deg V(h) \cap \bar{X} \leq D \deg X$.
As $V(h)\cap \bar{X}$ is a variety of pure dimension $s-1$ containing
$X$, the second claim follows by complete induction.

The last claim is obvious.

\section{Arithmetic varieties} \label{arvar}

Let $k$ be a number field with ring of integers $\CO_k$, and
$\Pe^t_{\CO_k}= \mbox{Proj}(\CO_k[x_0,\ldots,x_t])$ projective $t$-space over 
$\CO_k$.
A closed subset $\CX$ of $\Pe^t_{\CO_k}$ with induced scheme structure
will be called a subvariety if it is the
flat closure of a subvariety $X \subset \Pe^t_k$, i.e.\@ if $\CY$ is the
closure of a homogeneous ideal
$\fa \subset \CO_k[x_0,\ldots,x_t]$ with $\CO_k[x_0,\ldots,x_t]/\fa$ flat. 
A cycle with no finite part
is a linear combination with coefficients in $\Z$ of irreducible subvarieties. 
Write $X$ for the base extension of $\CX$ to
$\spec \; k$, and for any embedding $\sigma: k \hookrightarrow \C$ denote
by $X_{\sigma}$ the base extension of $X$ to $\C_\sigma$ as well as the
$\C$-valued points of $X$. If clear from the context, $X_\sigma$ will
also be denoted by $X$. We put $TR(\CX):=TR(X)$ the top regularity of $\CX$.

The canonical hermitian product on $k \otimes_{\Q} \C$ induces a metric on 
$O(1)$, and via arithmetic intersection theory for every cycle effective 
cycle $\CX$ on $\Pe^t_{\CO_k}$ the height $h(\CX) \in \R$ can be defined
(see e.g.\@ \cite{SABK}).

\satz{Proposition} \label{avar}
\begin{enumerate}
\item
\[ h(\Pe^t_{\CO_k}) = [k:\Q] \sigma_t, \quad \mbox{with} \quad 
   \sigma_t = \frac12 \sum_{k=1}^p \sum_{m=1}^k \frac1m \]
the $t$th Stoll number.

Further, for any effective cycle $\CY$ on $\Pe^t$ the height of
$\CY$ is nonnegative. Hence, if for two effective cycles $\CX,\CY$ in
$\Pe^t_{\CO_k}$ one denotes $\CX.\CY$ their intersection product, and
$X.Y$ the flat closure of the intersection product of their base
extension to $\spec \; k$, then
\[ h(X.Y) \leq h(\CX.\CY). \]

\item
Let $\CY \subset \Pe^t$ be a subvariety of codimension $p$ of $\Pe^t$,
and $f$ a global section of $O(D)$ whose restriction
to every irreducible component of $\CY$ is nonzero. Then, 
\[ h(\CY . \di f) = D h(\CY) + \int_{X(\C)} \log |f| \mu^{d-p} \delta_Y, \]
where $\mu$ is the first chern form of $\bar{O(1)}$, and the integral
is defined by resolution of singularities (See \cite{SABK}, II.1.2 for 
details).

\item
Let $\CX$ be a subvariety of codimension $p$ in projective space $\Pe^t$,
and $f \in \Gamma(\Pe^t,O(D))$ a global section. Then,
\begin{equation} \label{aanoben}
\int_X \log |f| \mu^{t-p} - \deg X \int_{\Pe^t} \log |f| \mu^t
\leq  c D \deg X, 
\end{equation}
with $c$ a positive constant only depending on $t$ and $p$.

\item
{\bf Arithmetic B\'zout Theorem:}
If $p,q$ are natural numbers with $p+q \leq t+1$, and
$\CX,\CY$ effective cycles in $\Pe^t$ of pure codimensions $p$ and $q$
intersecting properly, then
\[ h(\CX. \CY) \leq \deg Y h(\CX) + \deg X h(\CY) + \]
\[ [k:\Q] \left( \sigma_{2t+1-p-q} +\sigma_{t-p-q} + 
    \frac{p+q}2  \log 2 \right) \deg X \deg Y. \]

\item
Let $\CX \subset \Pe^t$ be an irreducible 
subvariety of codimension $p$, and \\ $f \in \Gamma(\Pe^t,O(D))_{\CO_k}$ 
a global section that has nonzero restriction to $\CX$, and
$f^\bot_X$ its projection to the orthogonal component of 
$I_\CX(D)$, the homogeneous polynomials of degree $D$ that are contained
in the prime ideal corresponding to $\CX$. Then, 

\begin{equation} \label{gsbezout}
h(\CX . \di f) \leq D h(\CX) + \deg X \log |f_X^\bot|_{L^2(\Pe^t)} + 
c D \deg X,
\end{equation}
where $c$ is a constant only depending on $t$, and the dimension of $X$.

\end{enumerate}
\end{Satz}

\proof
1. (\cite{BGS}, Proposition 3.2.4)

\vspace{2mm}

2.
\cite{BGS}, Proposition 3.2.1.iv or \cite{App1}, Proposition 3.8.1.

\vspace{2mm}

3. For $f$ with $\int_{\Pe^t} \log |f| \mu^t =0$, this is \cite{BGS},
5.1, remark (i). For general $f$ it also holds because the left hand
side of the inequality does not change if one replaces $f$ by
$af$ with $a \in \R$.

\vspace{2mm}

4. \cite{BGS}, Theorem 4.2.3.

\vspace{2mm}

5. Firstly, by part 2,
\[ h(\CX . \di f) = D h(\CX) + \int_X \log |f| \mu^{t-p}, \]
where $\mu$ is the Fubini-Study metric on $\Pe^t$, or alternatively the
first chern form of $\overline{O(1)}$.
Next, $f = f_X^\bot + g$ with $g \in I_X(D)_\C$. Hence,
\[ \int_X \log |f| \mu^{t-p} = \int \log |f_X^\bot| \mu^{t-p}, \]
which by (\ref{aanoben}) is less or equal
\[ \deg X \int_{\Pe^t} \log |f_X^\bot| \mu^t + c D \deg X \leq
   \frac12 \deg X \log \int_{\Pe^t} |f_X^\bot|^2 \mu^t + c D \deg X. \]

\subsection{Arithmetic bundles}

Let $\CO_k$ be a number ring. An arithmetic bundle over $\spec \; \CO_k$ is
a projective finitely generated $\CO_k$-module $E$ together with a 
hermitian product $\la \cdot |\cdot\ra$ on $E_\infty = E \otimes_{\CO_k} \C$.
For $\bar{F}$ a one dimensional arithmetic bundle define the arithmetic
degree
\[ \wde \; \bar{F} := \sum_{s \in S} - \log |v|_s, \]
where $v \in E$ is any nonzero element, and $S$ denotes the set
of all places of $\CO_k$. For an arbitrary arithmetic bundle define
\[ \wde \; \bar{E} := \wde \det \bar{F}. \]
If $\CO_k = \Z$, then $\wde \bar{E}$ is just minus the logarithm of
the covolume of $E$ in $E \otimes_\Z \R$.

If $\bar{E}$ is an arithmetic bundle, and $F \subset E$ a subbundle
the metric on $E_\infty$ induces a metric on $F_\infty$, and one
obtains an arithmetic bundle $\bar{F}$. If one uses the canonical isomorphism
of $(E/F)_\infty$ with the orthogonal complement $F^\bot_\infty$ of
$F_\infty$ in $E_\infty$ the metric on $E_\infty$
induces a metric on $(E/F)_\infty$, and one obtains an arithmetic bundle
$\overline{E/F}$. 
%Further, under suitable conditions on the metric
%on $E$, for example for the canonical metric on $\C^n=\Z^n \otimes_\Z \C$, 
%the module $G= E \cap F_\infty^\bot$ as rank $\rk E- \rk F$, and
%the metric on $E_\infty$ induces a metric on $G_\infty = F_\infty^\bot$. We
%have
%\begin{equation} \label{orth}
%\wde \bar{G} \leq \wde \; \overline{E/F}. 
%\end{equation}

\satz{Theorem(Minkowski)} \label{Minkowski}
Let $\bar{M}$ be an arithmetic bundle over $\spec \; \Z$, and
$K \subset M_{\otimes \Z} \R$ any closed convex subset that is symmetric
with respect to the origin, and fulfills
\[ \log vol(K) \geq - \wde(\bar{M}) + \rk M \log 2. \]
Then $K \cap M$ contains a nonzero vector. In particular taking
$K$ as the cube with logarithmic length of edge
$-\frac{1}{\rk M} \wde (\bar{M}) + \log 2$ centered at the origin, one sees
that there is a non zero lattice point $v \in M$ of logarithmic length
\[ \log |v| \leq - \frac1{\rk M} \; \wde (\bar{M}) + \frac12 \log \rk M. \]
\end{Satz}

\vspace{2mm}

Let now $\Pe^t_\Z = \Pe(\Z^{t+1})$ be projective space of dimension $t$, and
\[ E_D := \Gamma(\Pe^t,O(D)). \]
As $E_D = Sym^D E_1$, which in turn equals the space of homogeneous
polynomials of degree $D$ in $t+1$ variables, this lattice canonically 
carries the following metrics:

\begin{enumerate}

\item
The subspace metric $Sym^D E_1 \subset E_1^{\otimes D}$.

\item
The quotient metric $E_1^{\otimes D} \to Sym^D E_1$ referred to as $|f|_q$.

\item
The $L^2$-metric
\[ |f|^2_{L^2(\Pe^t)} = \int_{\Pe^t(\C)} |f|^2 \mu^t, \]
where $\mu$ is the Fubini-Study metric on $\Pe^t$.

\item
The supremum metric
\[ |f|_\infty = \sup_{x \in \Pe^t} |f_x|. \]

\end{enumerate}

Among these metric the relations 

\begin{equation} \label{Normrel}
 \log |f|_\infty - \frac D2 \sum_{m=1}^t \frac 1m \leq
   \int_{\Pe^t_\C} \log |f| \mu^t \leq \log |f|_{L^2} \leq
   \log |f|_\infty. 
\end{equation}
hold (\cite{BGS}, (1.4.10) or \cite{App1}, Lemma 3.1.).

\satz{Lemma} \label{basisnorm}
Let $I = (i_0, \dots i_t)$ be a multiindex with $|I| = i_0 + \cdots +i_t =D$, and
$X^I$ the polynomial $x_0^{i_0} \cdots x_t^{i_t}$. The set
$\{ X^I| \; |I| = D \}$ forms a basis of $\Gamma(\Pe^t_\Z,O(D))$, and
\[ |X^I|^2_{L^2} = |X^I|^2_q {D+t \choose t}^{-1} = 
   {D+t \choose I}^{-1} {D+t \choose t}^{-1}=
   \frac{i_0! \cdots i_t! t!}{(D+t)!} \leq 1. \]
Further, there is a positive constant $c_1$ depending on $t$ such that
\[ -c_1 D \leq \log |X^I|_{L^2} = 
   - \log {D+t \choose t} - \log {D+t \choose I} \leq 0, \]
and
\[ \sum_{\{ I| \; |I|=D\}} 
   \log |X^I|_{L^2} = - \sigma_t \frac{D^{t+1}}{(t+1)!} +
   O(D^t \log D). \]
\end{Satz}

\proof
The chain of equalities follows form \cite{BGS}, Lemma 4.3.6, and
its proof. The estimates for 
$-\log {D+t \choose t} - \log {D+t \choose I}$ are easy calculations.
The last equality is the arithmetic Hilbert-Samuel formula for $\Pe^t$.

\satz{Lemma} \label{polprod}
Let $f \in \Gamma(\Pe^t,O(D)), g \in \Gamma(\Pe^t,O(D'))$. Then,
\[ \log |f|_{L^2} + \log |g|_{L^2} -c_1(D+D') \leq
   \log |f g|_{L^2} \leq \]
\[ \log |f|_{L^2} + \log |g|_{L^2}+ c_1(D+D') \leq \]
\[ |f|_{L^2} + \log |g|_{L^2}+ c'_1(D+D')+
   \log {D+t \choose t}+ \log {D'+t \choose t}, \]
with $c_1'$ a positive constant only depending on $t$.
\end{Satz}

\proof
Let $I$ be any multiindex of order $D$. The vector
$X_I=x_0^{\otimes i_0} \otimes \cdots \otimes x_t^{\otimes i_t} \in 
E_1^{\otimes D}$ has length one, and by the previous Lemma,
\[ \log |X_I| = 0 \geq 
   \log |X^I|_q  \geq \log |X_I|- c_1 D . \]
Likewise, if $I'$ is a multiindex with $|I'|=D'$, then
\[ 0 \geq \log |X^{I'}|_q  \geq - c_1 D'. \]
Next, $|X_I \otimes X_{I'}| = |X_I| \; |X_{I'}|$, and thus
\[ \log |X^I X^{I'}|_q \leq \log |X_I X_{I'}| = 0 \leq 
   \log |X^I|_q + \log |X^{I'}|_q + c_1(D+D'), \]
hence, by the previous Lemma,
\[ \log |X^I X^{I'}|_{L^2} \leq \log |X^I|_{L^2} + \log |X^{I'}|_{L^2}+
   c_1(D+D') - \]
\[ \log {D+D'+t \choose t}+ \log {D+t  \choose t} + 
   \log{D'+t \choose t}. \]
Further, since $|X^I|_{L^2} \leq 1$,
\[ \log |X^I X^{I'}|_{L^2} \geq - c_1(D+D')\geq\log |X^I|_{L^2} +
   \log |X^{I'}|_{L^2} - c_1 (D+D'), \]
proving the claim for $f$, and $g$ monomials. The claim follows for
general $f$, and $g$ because the $X^I,X_I$ form orthogonal bases in 
any of the above metrics.

\vspace{2mm}

Let $X$ be a subvariety of pure dimension $s$ in $\Pe^t_\C$.
Then on 
\[ I_X(D) := \{ f \in H^0(\Pe^t,O(D)) | \; |f|_X = 0 \}, \]
there are the restrictions of the norms $|\cdot|_{Sym}, |\cdot|_q$, and
$|\cdot|_{L^2(\Pe^t)}$, and on
\begin{equation} \label{quotient}
F_\CX(D) = H^0(\CX,O(D)), 
\end{equation}
there is the quotient norm $|\cdot|_q$ induced by the quotient norm
$|\cdot|_q$ via the canonical quotient map
\[ q_D: E_D \to F_X(D), \]
the $L^2(\Pe^t)$-norm
\[ | \cdot |_{L^2(\Pe^t)}: F_D(X) \to \R,  \quad
         \bar{f} \mapsto \inf_{q_D(f)=\bar{f}} 
         \sqrt{\int_{\Pe^t} |f|^2 \mu^t} =
         \sqrt{\int_{\Pe^t} |f_X^\bot|^2 \mu^t} = |f_X^\bot|_{L^2(\Pe^t)}, \]
with $f_X^\bot$ the unique vector that is orthogonal to
$I_X(D)$, and fulfills $q_D(f_X^\bot) = \bar{f}$,
and the $L^2(X)$-norm
\[ |\cdot|_{L^2(X)}: F_D(X) \to \R, \quad f \mapsto \int_X |f| \mu^p. \]

\satz{Lemma} \label{polprodqu}
Let $X\subset \Pe^t_\C$ be a subvariety, and
$\bar{f}\in F_D(X),\bar{g} \in F_{D'}(X)$. Then,
\[ \log |\overline{f g}|_{L^2(\Pe^t)} \leq 
   \log |\bar{f}|_{L^2(\Pe^t)} + \log |\bar{g}|_{L^2(\Pe^t)} + c_1(D+D') \leq \]
\[ \log |\bar{f}|_{L^2(\Pe^t)} + \log |\bar{g}|_{L^2(\Pe^t)} +
   c_1'(D+D')+ \log {D+t \choose t} + \log {D'+t \choose t}. \]
\end{Satz}

\proof
Let $f\in \Gamma(\Pe^t,O(D), g \in \Gamma(\Pe^t,O(D'))$ be representatives
of $\bar{f},\bar{g}$, and
$f=f_1+f_2$ with  $f_1 \in I_X(D)$, and $f_2 \in I_X(D)^\bot$, and
likewise for $g$. Then, $|\bar{f}|_{L^2(\Pe^t)} = |f_2|_{L^2},
|\bar{g}|_{L^2(\Pe^t)} = |g_2|_{L^2}$. Let further
$\overline{f g}$ be represented by $h = fg \in \Gamma(\Pe^t,O(D+D'))$ 
with decomposition $h=h_1 + h_2$. Then,
$h_2 = h-h_1 = (f_1 g_1 + f_1g_2 + f_2 g_1 -h_1) + f_2 g_2$
with $f_1 g_1 + f_1g_2 + f_2 g_1 -h_1 \in I_X(D+D')$, and 
$h_2 \in I_X(D+D')^\bot$. Further, $\overline{f_2 g_2} = \overline{fg}$.
Consequently, by the previous Lemma,
\[ \log |\overline{f g}|_{L^2(\Pe^t)} = \log |h_2|_{L^2} \leq
   \log |f_2g_2|_{L^2} \leq \]
\[ \log |f_2|_{L^2} + \log |g_2|_{L^2} + 
   c_1(D+D') +  \log {D+t \choose t} + \log {D'+t \choose t} = \] 
\[ \log |\bar{f}|_{L^2(\Pe^t)} + \log |\bar{g}|_{L^2(\Pe^t)} + c_1(D+D') + 
    \log {D+t \choose t} + \log {D'+t \choose t}. \]

\vspace{2mm}

\section{Arithmetic Hilbert functions}

For $\CX \subset \Pe^t-\Z$ a subvariety, let 
$F_\CX(D)=\Gamma(\CX,O(D))/I_\CY(D)$.
%and 
%$G_\CX(D) = \Gamma(\Pe^t,O(D))_\Z \cap I_X(D)_\infty$
%as in setction \ref{arvar}, number 5.
If not stated otherwise, in this section $|\cdot|$ will alway denote
the $L^2$-metric $|\cdot|_{L^2(\Pe^t)}$, on 
\[ E(D) = \Gamma(\Pe^t,O(D))_\R = \Gamma(\Pe^t,O(D))_\Z \otimes_\Z \R, \]
\[ I_X(D)\subset E(D), \quad F_X(D) = F_\CX(D) \otimes_\Z \R  
%\quad G_X(D)
. \]

For $\CX$ a subvariety of pure dimension $s+1$ of $\Pe^t_\Z$, define
the arithmetic Hilbert functions. 
\[ \hat{H}_\CX(D) := \wde (\bar{F}_\CX(D),|\cdot|_{L^2(\Pe^t)}), \quad
   \hat{\CH}_\CX(D) := \wde (\bar{F}_\CX(D),|\cdot|_{L^2(X)}). \]

\satz{Lemma} \label{mlnull}
Let $\bar{M} \subset\bar{N} \subset \bar{E}(D)=\overline{\Gamma(\Pe^t,O(D))}$ 
be arithmetic subbundles. Then, 
\[ \wde \; \overline{E(D)/M} \geq  \wde \; \overline{N/M} \geq
   \wde \; \bar{N}; \]
in particular
\[ \hat{H}_\CY(D) \geq 0 \]
for every subvariety $\CY$ of $\Pe^t$.
Also, with $d \leq \deg Y$ the number of irreducible components of $Y$,
\[ \hat{H}_{\CY}(D_1) \geq \hat{H}_\CY(D) - 
   (2 (D_1-D)\log d + c_1'D_1)H_Y(D), \]
for every $D_1 \geq D$.
\end{Satz}

\proof
Let $q:E(D) \to E(D)/M$ be the canonical projection.
Since $B= \{ X^I| \; |I| = D\}$ is a basis of $E(D)$, there is
as subset $\bar{B} \subset B$ such that $q(\bar{B})$ is as basis
of $E(D)/M$, and because of $|q(X^I)| \leq |X^I| \leq 1$, the first 
two inequalities and $\hat{H}_\CY(D) \geq 0$ follow.

Since for every irreducible component $\CY_j$ of $\CY$, there is an 
$i \in \{0,\ldots,t \}$ such that $x_i$ restricted to $\CY_j$ is nonzero,
by Lemma \ref{summe}, there is an $L \in \Gamma(\Pe^t, O(1))$, such
that $\log |L| \leq 2 \log d$, and the restriction of $L$ to
every irreducible component is nonzero. Hence, by Lemma \ref{polprodqu},
\[ \wde \;(L^{D_1-D} \bar{E}_\CY(D)) \geq \wde \bar{E}_\CY(D) - 
   (c_1' (D_1)+2(D_1-D) \log d) \rk \; E(D) = \]
\[ \hat{H}_\CY(D) - H_Y(D) (2(D_1-D)\log d +c_1' D_1).  \]
As $L^{D_1-D} \bar{E}_\CY(D)$ is a submodule of $\bar{E}_\CY(D_1)$, the last claim
follows.

\satz{Theorem} \label{arHilbert}
Let $\CX$ be a subvariety of pure dimension $s+1$ of $\Pe^t_\Z$.

\begin{enumerate}

\item
For every $D \in \N$,
\[ \hat{\CH}_\CX(D)\leq \deg X
   \left( D h(\CX) + \frac12 (\log \deg X+2s \log D)\right) {D+s \choose s}. \]

\item
With a positive constant $c_2'$ only depending on $t$, and $p$,
\[ \hat{H}_\CX(D) \leq \]
\[ \left( D h(\CX) + c_2' \deg X D +
   \deg X \left(\frac12 \log \deg X + s \log D\right)\right) {D+s \choose D}. \]
Hence for $X$ with $\deg X \leq D^t$, and with $c_2 = c_2'+2t$, $\deg X$
\[ \hat{H}_\CX(D) \leq (h(\CX)+c_2 \deg X) D {D+s \choose D} . \]

%\item
%With a positive constant $c'$ only depending on $t$, and $p$,
%\[ \hat{H}_\CX(D) \geq \wde (\bar{G}_\CX(D)) \geq
%   - \hat{H}_\CX(D) + \hat{H}_{\Pe^t}(D) \geq \hat{H}_{\Pe^t}(D) - \]
%\[ \left( D h(\CX) + 
%   \deg X D (\sigma_{t-1}-\sigma_t) +
%   \deg X (\frac12 \log \deg X + s \log D) \right) {D+s \choose D}. \]
%Hence for $X$ with $\deg X \leq D^t$, and $c=c'+t$,
%\[ \hat{H}_\CX(D) \geq -(h(\CX)+c_5\deg X) D {D+s \choose s}. \]

%\item
%There is a constant $c_8>0$ such that for any $(D,aD)$-regular 
%subvariety $\CX \subset \Pe^t$ of relative dimension $0$,
%\[ \hat{H}_\CX(D) \geq h(\CX) D - \deg X D - c_8 \log D. \]

\item
There are positive constants $c_3,c_3' > 0$ only depending on $t$ and $s$
such that if $\CY=\CX$ is irreducible, 
then for $D \geq 2 \; TR(X)$, the inequality
\[ \hat{H}_\CY(D) \geq (c_3 h(\CY) - c_3' \deg Y) D^{s+1} \]
holds.

\end{enumerate}
\end{Satz}

\proof
1. Let $\CX = \CX_1 \cup \cdots \cup \CX_n$ be the decomposition into
irreducible components. We use complete induction on $n$. If $n=1$, i.\@
e.\@ $\CX$ is irreducible, let $f \in F_\CX(D)$ be nonzero. 
Then, by Proposition \ref{avar}.2 and and the concavity of the log-function,
\begin{eqnarray*}
0 \leq h(\CX . \di (f)) &=& D h(\CX) + \int_X \log |f| \mu^s \leq
   D h(\CX) + \log \int_X |f| \mu^s \\
 &=& D h(\CX) + \log |f|_{L^2(X)}. 
\end{eqnarray*}
Hence,
\[ \log |f|_{L^2(X)} \geq - D h(\CX) \]
for every nonzero vector $f \in F_D(\CX)$.
By the Theorem of Minkowski, thus
\[ - \wde (\bar{F}_D,|\cdot|_{L^2(X)}) \geq -\mbox{rk} F_D D h(\CX) -
     \frac{rk F_D}2 \log rk F_D. \]
By Proposition \ref{alHilbert}.1, 
$- rk F_D \geq - \deg X {D+s \choose s}$; hence the above is greater or equal
\[ -\deg X
   \left( D h(\CX) + \frac12 (\log \deg X+ 2s\log D)\right) {D+s \choose s}, \]
which proves the claim for $n=1$.
Assume now the claim has been proved for $n-1$. We have the 
surjective restriction map
\[ \varphi: F_{\CX}(D) \to F_{\CX_1 \cup \cdots \cup \CX_{n-1}}(D), \]
and $\bar{F}_{\CX_1 \cup \cdots \cup \CX_{n-1}}(D) = \overline{F_{\CX}/
\mbox{ker} \varphi}$. Hence,
\[ \hat{\CH}_\CX(D) = \wde \bar{F}_{\CX}(D) =
   \wde \; \overline{\mbox{ker} \varphi} + 
   \wde \bar{F}_{\CX_1 \cup \cdots \cup \CX_{n-1}}(D), \]
which by induction hypothesis is at most
\[ \wde\;  \overline{\mbox{ker} \varphi} + 
   \sum_{i=1}^{n-1} \deg X_i
   \left( D \sum_{i=1}^{n-1} h(\CX_i) + 
   \frac12 \left(\log \left(\sum_{i=1}^{n-1} \deg X_i \right)+ 
   2s\log D\right)\right) {D+s \choose s}. \]
Since $\mbox{ker} \varphi$ maps injectively to $F_{\CX_n}(D)$, 
for every $f \in \mbox{ker f}$
\[ \log |f| \geq \log |f_{X_n}^\bot| \geq - D h(\CX_n), \]
and
\[ \rk \; \mbox{ker} \varphi \leq H_{X_n}(D) \leq \deg X_n {D + s \choose s}, \]
which together with the Theorem of Minkowski implies
\[ \wde \;\overline{\mbox{ker}\varphi} \leq \deg X_n \left(D h(\CX_n) +\frac12 
   (\log \deg X_n+ 2s\log D)\right) {D+s \choose s}. \]
Hence,
\[ \hat{\CH}(D) \leq \deg X_n \left( D h(\CX_n) + \frac12 
   (\log \deg X_n+ 2s\log D)\right) {D+s \choose s} + \]
\[  \sum_{i=1}^{n-1} \deg X_i
   \left( D \sum_{i=1}^{n-1} h(\CX_i) + 
   \frac12 \left(\log \left(\sum_{i=1}^{n-1} \deg X_i \right)+ 
   2s\log D\right)\right) {D+s \choose s} \leq \]
\[ \deg X
   \left( D h(\CX) + \frac12 (\log \deg X+2s \log D)\right) {D+s \choose s}, \]
finishing the proof.

\vspace{2mm}

2.
Let $\CX = \CX_1 \cup \cdots \cup \CX_n$ be the decomposition into
irreducible components. If $n=1$, i.\@ e.\@ $\CX$ is irreducible, denote by  
$f_X^\bot$ the orthogonal projection of $f \in \Gamma(\Pe^t,O(D))$
modulo $I_{X_i}(D)$. Clearly $|f_X^\bot| \leq |f|$.
Further by Proposition \ref{avar}.5,
\[ \deg X \log |f_X^\bot|_{L_2(\Pe^t)} \geq 
   h(\CX . \di(f)) - D h(\CX) - c D \deg X \geq -D h(\CX) - c D \deg X. \]
Consequently,
\[ \log |f|_{L^2(\Pe^t)} \geq - D \frac{h(\CX)}{\deg X} - c D . \]
As $f \in I_D(\CX)^\bot \setminus \{0 \}$ was arbitrary, the claim
follows for $n=1$ in the same way as part one.

Assume now the claim has been proved for $n-1$.
With the notations of the proof of part one, 
\[ \hat{H}_\CX(D) = \wde \; \overline{\mbox{ker}\varphi} +
   \hat{H}_{\CX_1 \cup \cdots \cup \CX_{n-1}}(D), \]
which by induction hypothesis is less or equal 
\[ \wde \; \overline{\mbox{ker}\varphi} + \]
\[ \left( D \sum_{i=1}^{n-1} h(\CX_i) + c_2' \sum_{i=1}^{n-1} \deg X_i D +
   \sum_{i=1}^{n-1} \deg X_i \left(\frac12 \log 
   \sum_{i=1}^{n-1} \deg X_i + s \log D\right) \right) {D+s \choose D}. \]
Similarly as in the proof of part 1, one proves
\[ \wde \; \overline{\mbox{ker} \varphi} \leq
    \left( D h(\CX_n) + c_2' \deg X_n D +
   \deg X_n \left(\frac12 \log \deg X_n + s \log D\right) \right) 
   {D+s \choose D}, \]
and part 2 thus follows by taking the sum of the two terms.

\vspace{2mm}

\subsection{Auxiliary varieties}

The proof of the last part of Theorem \ref{arHilbert} requires some 
preparations.

%\vspace{2mm}

%\satz{Corollary} \label{isoliert}
%Let $\CY \subset \Pe^t$ be an irreducible subvariety of codimension $s$, and 
%$\CY_1,\ldots,\CY_k$ irreducible varieties of codimension $r<s$ that contain
%$\CY$ and are each a locally complete intersection of hypersurfaces
%of degree at most $D$. Assume further that
%for some $l \in \N$ every $i,j =1,\ldots,k$ the the inequality
%$\deg Y_i \leq \deg Y_j$ holds and the 
%the inclusion $I_{Y_i}(lD-1) \hookrightarrow I_Y(lD-1)$ is an isomorphism
%for every $i=1,\ldots,k$. 

%Then $\CY_i = \CY_j$ for every $i,j=1,\ldots$ and
%\[ h(\CY_i) \leq c \frac{h(\CY) + c_5 \deg Y}{(lD)^{s-r}}, \]
%where $c$ is a constant only depending on $r,s,t$, but not on $l$ and $D$.
%\end{Satz}

%\proof

\vspace{2mm}

\satz{Lemma} \label{abschn}
Let $\CY \subset \Pe^t_{\Z}$ be an irreducible variety of codimension $s$, 
and $\bar{D}= TR(Y)$ the top regularity of $Y$. For $m \geq 2t$,
assume $\bar{D} \geq 2 m^s$,
set $D_r := \left[\frac{\bar{D}}{m^{s+1-r}}\right]$, and let $r \geq 0$ be 
the maximal number such that there is a subvariety $X_r$ of pure
codimension $r$ that is a locally complete intersection of 
hypersurfaces of degree at most $D_r$ at $\CY$. Then, there is an
irreducible component $\bar{\CY}=\CY_r$ of $\CX_r$ such that
$I_\CY(D_{r+1}) = I_{\bar{\CY}}(D_{r+1})$ and for every irreducible component
$Z$ of $X_r$,
\[ \frac 1{(2t)^{2t} t! m^{(s-r)(s-r-1)}}\bar{D}^{s-r} \deg \bar{Y} \leq \deg Y \leq 
   \bar{D}^{s-r} \deg Z. \]
Also, $r<s$ if $s \geq 1$.
%\item
%if Theorem \ref{arHilbert}.4 holds for varieties of codimension at most
%$s-1$, and $m\geq N$ in Theorem \ref{arHilbert}, then also
%\[ h(\bar{\CY}) \leq (c' h(\CY) + c'' \deg Y) D^{r-s}, \]
%where $c,c',c''$ are positive constants only depending on $r,s,t$ and $m$. 

\end{Satz}

\vspace{2mm}

{\bf Remark:} For $m$ bigger than a given constant only depending on $t$,
it can be achieved
that $\CX_r$ is irreducible, hence equal to $\bar{\CY}$, but this fact
will not be needed in the following.

\proof
%
%I claim $\deg Z \leq c \deg Y D^{r-s}$ for every irreducible component $\CZ$ of 
%$\CX_r$. Let $M \leq D_r$ be the smallest number such that $\CZ$ is an 
%irreducible component of a locally complete intersection of degree at most 
%$M$, and $M_k = \left[\frac M{m^k}\right]$. Assume first 
%$m M \geq D_r$. As $\mbox{codim} \CZ \leq s-1$, with
%$r'<r$ the biggest number such that there is a subvariety 
%$\CX_{r'} = \CZ_1 \cup \cdots \cup\CZ_l$ of
%pure comdimension $r'$ that is a locally complete intersection of 
%hypersurfaces of degree at most $M_k$, the induction hypothesis implies
%$\deg Z_j \geq \deg Z M^{k-r}$ for all $j \in \{1,\ldots,l\}$.
%
%By definition of $M_{k+1}-1$ and Lemma \ref{summe}, 
%$I_Z(M_{k+1}-1)\setminus I_{Z_j}(M_{k+1}-1) = \emptyset$ for some 
%$j \in \{1,\ldots,l\}$. Since
%$I_Z \subset I_Y$ also $I_Y(M_{k+1}-1) \setminus I_{Z_j}(M_{j+1}) = \emptyset$,
%that is
%\[ 0 = \dim I_Y(M_{k+1}-1)/ I_{Z_j}(M_{k+1}) =
%   H_{Z_j}(M_{k+1}-1) - H_Y(M_{k+1}-1) \geq \]
%\[ \deg Z_j {M_{k+1}-1 -k M_k +t+k \choose t-k} - 
%   \deg Y {M_{k+1}-1 + t-s \choose t-s} \geq \]
%\[ \deg Z_j \left(1-\frac{k}m\right)^{t-k} 
%   \frac{M_{k+1}^{t-k}}{(t-k)!}
%   - \deg Y (2(t-s))^{t-s} \frac{M_{k+1}^{t-s}}{(t-s)!}. \]
%For $m \geq 2t$ this implies
%\[ \deg Y \geq \frac1{2^{t-k}(2(t-s))^{t-s}} \frac{(t-s)!}{(t-k)!}
%   \deg Z_j M_{k+1}^{s-k} \geq \]
%\[ \frac1{2^{t-k}(2(t-s))^{t-s}} \frac{(t-s)!}{(t-k)!} \deg Y M_{k+1}^{s-r}=
%   \frac1{2^{t-k}(2(t-s))^{t-s}} \frac{(t-s)!}{(t-k)!} 
%   \frac{\deg Y}{m^{(s-r)(s+1-r)}} D^{s-r}, \]
%implying the claim with $c=\cdots$ suitably chosen.
%
 Let $\CX_r=\CZ_1 \cup \cdots \cup \CZ_l$ be the decomposition into 
irreducible varieties and assume
$I_{\CY}(D_{r+1}) \neq I_{\CZ_j}(D_{r+1})$ for all $j \in \{1,\ldots,l\}$. Choose
an $h_j \in I_{\CY}(D_{r+1})\setminus I_{\CZ_j}(D_{r+1})$ for each $j$. 
Then, by Corollary \ref{summe1}, there
is a polynomial $h$ of degree $D_{r+1}$, such that $\di\;\! h$ intersects
$\CX_r$ properly. Hence, by Lemma \ref{bas}.3 and Corollary \ref{TR}, 
there is a variety
$\CX_{r+1}$ that is a locally complete intersection of homogeneous polynomials
of degree at most $D_{r+1}$ at $\CY$, contradicting the maximality of $r$.

By a similar argument, it follows $r<s$ if $s\geq 1$ since
$D_r = \left[ \frac{\bar{D}}{m^{s+1-r}} \right] <\bar{D}$.

Choose an irreducible component $\bar{\CY}$ of $\CX_r$ such that
$I_\CY(D_{r+1}) = I_{\bar{\CY}}(D_{r+1})$. Then 
$TR(\bar{Y}) \leq D_r$, hence by Theorem \ref{alHilbert},
\begin{eqnarray*} 
 0 &=& \dim I_Y(D_{r+1})/ I_{\bar{Y}}(D_{r+1}) = H_{\bar{Y}}(D_{r+1}) - H_Y(D_{r+1}) \\
   &\geq& \deg \bar{Y} {D_{r+1} -r D_r +t-r \choose t-r} - 
   \deg Y {D_{r+1} + t-s \choose t-s} \\
   &\geq& \deg \bar{Y} \left(1-\frac{r}m\right)^{t-r} 
   \frac{D_{r+1}^{t-r}}{(t-r)!}
   - \deg Y (2(t-s))^{t-s} \frac{D_{r+1}^{t-s}}{(t-s)!}. 
\end{eqnarray*}
As $m \geq 2t$, this implies
\[ \deg Y \geq \frac1{2^{t-r}(2(t-s))^{t-s}} \frac{(t-s)!}{(t-r)!}
   \deg \bar{Y} D_{r+1}^{s-r} \geq \]
\[ \frac1{2^{t+s-2r}(2(t-s))^{t-s}} \frac{(t-s)!}{(t-r)!}
   \frac{\deg \bar{Y}}{m^{(s-r)(s-r-1)}} \bar{D}^{s-r}. \]
Further, by Corollary \ref{TR}
\[ \deg Y \leq \deg Z \bar{D}^{s-r}, \]
finishing the proof.

%\vspace{2mm}

%2. By part 1,
%\[ \hat{H}_{\bar{\CY}}(D_{r+1}) = 
%   \wde(\overline{\Gamma(\Pe^t),O(D_{r+1})})- \wde(\bar{I}_{\bar{CY}}(D_{r+1})) =
%   \wde(\overline{\Gamma(\Pe^t),O(D_{r+1})})- \wde(\bar{I}_{\CY}(D_{r+1})) =
%   \hat{H}_{\CY}(D_{r+1}). \]
%If $s=0$, nothing is to be proved. If $s\geq 1$, then $r\leq s-1$, hence
%Theorem \ref{arHilbert}.4 can be applied to $\bar{\CY}$, and of course
%Theorem \ref{arHilbert}.2 can be applied to $\CY$ leading
%\[ (c_6(h(\bar{\CY}) - c_7 \deg X) D_{r+1} {D_{r+1}+t-r \choose t-r} \leq
%   h(\CY) + c_5 \deg X) D_{r+1} {D_{r+1}+t-s \choose t-s}. \]
%Since $D_{r+1}=\la \frac{D}{m^{r+1-s}} \ra$, the claim follows with suitable
%choices for $c',c''$.

\satz{Proposition} \label{chain}
Let $\CY \subset \Pe^t$ be an irreducible subvariety of codimension $s$, and
$\bar{D} =TR(\CY)$. Let further $r$,
$\bar{\CY}$ be as in the previous Lemma, $n \in \N$, 
and $D = \bar{D}+D_{r+1}$. There
are constants $c_4(k),c_4'(k)$ only depending on $k,t,r,s$ such that if 
\[ \hat{H}_{\CY}(D) \leq \frac1n (h(\CY) + \deg Y) D^{t+1-s}, \]
then there is a chain of irreducible subvarieties
\[ \CY = \CY_s \subset \CY_{s-1} \subset \cdots \subset \CY_{r+1} \subset \CY_r = 
   \bar{\CY}, \]
such that  for $k \leq s-1$,
\[ \deg Y \leq \deg Y_k \bar{D}^{s-k} \leq \deg \bar{Y} (2\bar{D})^{k-r}, 
   \quad \mbox{and} \]
\[ h(\CY_k) \leq \frac {c_4(k)}n\frac{h(\CY)+\deg Y}{\bar{D}^{s-k}} +
   c_4'(k) h(\bar{\CY}) \bar{D}^{k-r}. \]
\end{Satz}

\proof
For $k=r$, take $\CY_r = \bar{\CY}$, $c_4'(r)=1$, and $c_4(r)=0$. 
For $k\leq s-2$, 
assume $\CY_k$ with the required property
has been construed. Firstly, since $Y_k \supset Y$, one has
$I_{Y_k}(D_{r+1}) \subset I_Y(D_{r+1})$ which by the Lemma equals
$I_{\bar{Y}}(D_{r+1})$. Hence, $H_{Y_k}(D_{r+1}) = H_{\bar{Y}}(D_{r+1})$.
Since $TR(\bar{Y}) \leq D_r < \frac{D_{r+1}}{m-1}$, Proposition
\ref{alHilbert} implies
\[ H_{Y_k}(D_{r+1}) = H_{\bar{Y}}(D_{r+1}) \geq 
   \deg \bar{Y} \frac{D_{r+1}^{t-r}}{2^{t-r}(t-r)!} \geq  \]
\[ \deg Y_k  \frac1{2^{t+k-2r}m^{(t-r)(s-r)}(t-r)!} \bar{D}^{t-k}, \]
the last inequality following from the induction hypothesis.
Next, since \\
$Y = \overline{U \cap V(h_1) \cap \cdots \cap V(h_s)}$ with
$\deg h_j \leq \bar{D}$ for  $j=1,\ldots,s$ 
there is an $i \in \{1,\ldots,s\}$ such that
$h_i \in I_Y(\bar{D}) \setminus I_{Y_k}(\bar{D})$, hence the map
\[ \Gamma(Y_{k+1},O(D_{r+1})) \to I_Y(D)/I_{Y_k}(D), \quad
   f \mapsto f \cdot h_i \]
is an injection, implying
\[ \rk \; I_Y(D)/I_{Y_k}(D) \geq  H_{Y_k}(D_{r+1}) \geq \deg Y_k 
   \frac1{2^{t+k-2r}m^{(t-r)(s-r)}(t-r)-!} D^{t-k}. \]
Further, by assumption, and Lemma \ref{mlnull},
\[ -\wde \overline{I_{\CY}(D)/I_{\CY_k}(D)} =
   \hat{H}_{\CY}(D) - \hat{H}_{\CY_k}(D) \leq
   \frac1n (h(\CY) + \deg Y) D^{t+1-s}. \]
Finally,
\[ \rk \;I_Y(D)/I_{Y_k}(D) \leq H_{\Pe^t}(2\bar{D}) = {2\bar{D}+t \choose t}, \]
hence
\[ \log \rk \; I_Y(D)/I_{Y_k}(D) \leq t \log(2\bar{D}+t) \leq 
   t\log(3t\bar{D}) \]
   
Thus, by the Theorem of Minkowski, there is an 
$f \in I_{\CY}(D)/I_{\CY_k}(D)$ such that
\begin{eqnarray*}
\log |f_{\CY_k}^\bot| &\leq &
   \frac{t!2^{t+k-2r}m^{(t-r)(s-r)}}{n} 
   \frac{h(\CY)+\deg Y}{\deg Y_k} D^{k+1-s} + \frac t2 \log(3tD) \\
  &\leq& \frac{t!2^{2t+k+1-3r}m^{(t-r)(s-r)}}{n} 
   \frac{h(\CY)+\deg Y}{\deg Y_k} D^{k+1-s}. 
\end{eqnarray*}
Choose $\CY_{k+1}$ as any irreducible component of $\CY_k . \di f$ that
contains $\CY$. Then, by the Theorem of B\'ezout and induction hypothesis,
\[ \deg Y_{k+1} \leq \deg Y_k \deg f = \deg \bar{Y} D(2\bar{D})^{k-r} \leq
   \deg \bar{Y} (2\bar{D})^{k+1-r}. \]
Also, by Proposition \ref{avar}.5 and induction hypothesis,
\begin{eqnarray*}
  h(\CY_{k+1}) &\leq&  2\bar{D} h(\CY_k) + \deg Y_k \log |f_{\CY_k}^\bot| +
   2c\bar{D} \deg Y_k \\
  &\leq& \frac{2c_4(k)}n \frac{h(\CY) + \deg Y}{\bar{D}^{s-(k+1)}} +
   2c_4'(k)h(\bar{\CY}) \bar{D}^{k+1-r} \\
 &+&  \frac{t!2^{2t+2k+1-3r}m^{(t-r)(s-r)}}{n!} 
   (h(\CY)+\deg Y) \bar{D}^{k+1-s} + c \deg \bar{Y}(2\bar{D})^{k+1-r}. 
\end{eqnarray*}
Since $\deg \bar{Y} \leq (2t)^{2t} t!m^{(s-r)(s-r-1)}\deg Y\bar{D}^{r-s}$, 
for $m=2t$ and suitably chosen $c_4(k+1),c_4'(k+1)$ the above is less or equal
\[ \frac{c_4(k+1)}n (h(\CY)+\deg Y) \bar{D}^{k+1-s} + 
   c'_4(k+1) h(\bar{\CY}) \bar{D}^{k+1-r}, \]
finishing the proof.

\vspace{2mm}

Put now 
\[ m=2t, \quad e_0= \frac 1m,  \quad e_1=1+\frac1m,\quad
   e_2=\frac1{(4t)^{t}t!}, \quad
   e_3 = \frac1{(2t)^{2t} t! m^{(s-r)(s-r-1)}}, \]
\[ e_4=\frac1 {(4t)^{2t} t! m^{(s-r)(s-r-1)}}, \quad e_5 = c_4'(s-1)  \]
and $n=8(\left[c_4(s-1)\right]+1)$.

\satz{Corollary} \label{Ysm1}
With the above constants, let $\CY \subset \Pe^t$ be an irreducible variety of
codimension $s \geq 1$ and $\bar{D}= TR(\CY)$, and
assume that with $D_0 = [e_0\bar{D}],\tilde{D} = [e_1 \bar{D}]$,
\[ \hat{H}_{\CY}(\tilde{D}) \leq e_2 (h(\CY) - \deg Y) \tilde{D}^{t+1-s}. \]
Then, there exist irreducible subvarieties $\bar{\CY},\CY_{s-1}$ with the
following properties:
\begin{enumerate}
\item
\[ \rk \; I_Y(\tilde{D})/I_{Y_{s-1}}(\tilde{D}) \geq e_2 \deg Y \tilde{D}^{t-s}, \]
and 
\[ I_Y(D_0) = I_{\bar{Y}}(D_0) = I_{Y_{s-1}}(D_0). \]
\item
$\bar{\CY}$ is an irreducible component of a locally complete intersection
of degree at most $\frac{\tilde{D}}{2t}$ at $\CY$ and codimension $r<s$. Further,
\[ e_3 \deg \bar{Y} \bar{D}^{s-r} \leq \deg Y \leq \deg \bar{Y} \bar{D}^{s-r}, \]
\item
$\CY_{s-1}$ has codimension $s-1$, contains $\CY$ and fulfills
\[ e_4 \deg Y \leq e_4 \bar{D} \deg Y_{s-1} \leq 2 \deg Y, \]
and
\[ h(\CY_{s-1}) \leq \frac18 \frac{h(\CY)}{\bar{D}} + 
   e_5 h(\bar{\CY}) \bar{D}^{s-1-r} +
   \frac18 \frac{\deg Y}{\bar{D}}. \] 
\end{enumerate}
\end{Satz}

\proof
Let $\bar{Y}$ be the variety from Lemma \ref{abschn} and $Y_{s-1}$ as in 
Proposition \ref{chain}.

1. By Lemma \ref{abschn},
$I_{Y{s-1}}(D_0) \subset I_Y(D_0) = I_{\bar{Y}}(D_0) \subset I_{Y_{s-1}}(D_0)$.
Hence, just as in the proof of Proposition \ref{chain}, we have
\[ \rk I_Y(\tilde{D})/I_{Y_{s-1}}(\tilde{D}) \geq H_{Y_{s-1}}(D_0) = 
   H_{\bar{Y}}(D_0) \geq
   \frac 1{2^t t!}\deg \bar{Y} D_0^{t-r} \geq \]
\[ \frac 1{(4t)^{t}t!} \deg \bar{Y} \bar{D}^{t-r} \geq e_2 \deg Y \bar{D}^{t-s}. \]

\vspace{2mm} 
  
2. Is a reformulation of inequalities in Lemma \ref{abschn}.

\vspace{2mm}

3. Is a reformulation of the inequalities in Proposition \ref{chain}
for $k=s-1$, and Lemma \ref{abschn}.

\subsection{The lower bound}

\satz{Lemma} \label{mldrei}
Let $\CY \subset \CX \subset \Pe^t$ be irreducible
subvarieties of codimensions $s$, and $s-1$ respectively. 
For every $f \in I_\CY(D) \setminus I_{\CX}(D)$, let again
$f_X^\bot$ be the orthogonal projection of $f$ modulo $I_X(D)$. Then,
with some positive constant $c_5$ only depending on $s,t$,
\[ \log |f_X^\bot|_{L^2(X)} = \log |f|_{L^2(X)} \geq h(\CY) - D h(\CX) 
   \quad \mbox{and} \]
\[ \log |f|_{L^2(\Pe^t)} \geq \log |f_X^\bot|_{L^2(\Pe^t)}
   \geq \frac{h(\CY)}{\deg X} - 
   D \frac{h(\CX)}{\deg X} - c_5 D \]
for every $f \in I_\CY(D)$ with $f_X^\bot \neq 0$.
\end{Satz}

\proof 1.
Let $\CZ = \CX . \di f$. Then, clearly $\CY$ is a subvariety of codimension $0$
of $\CZ$, hence $h(\CY) \leq h(\CZ)$.
By Proposition \ref{avar}.2, and the concavity of the $\log$-function,
\[ h(\CY) \leq h(\CZ) =
   D h(\CX) + \int_X \log |f| \mu^{t+1-s} \leq 
   D h(\CX) + \log |f|_{L^2(X))}, \]
and the first claim follows.
Further, by Proposition \ref{avar}.3 with $c_5$ from (\ref{aanoben}),
\[ \int_X \log |f| \mu^{t+1-s} = \int_X \log |f_X^\bot| \mu^{t+1-s} \leq  
   \deg X \int_{\Pe^t} \log |f_X^\bot| \mu^t + c_5 D \deg X \leq \]
\[ c_5 \deg X D + \deg X \log |f_X^\bot|_{L^2(\Pe^t)}, \]
proving the second claim.

%\satz{Lemma} \label{mlvier}
%Let $\CY$ be an $(\bar{D},a\bar{D})$ regular, irreducible subvariety of 
%codimension $s$ of $\Pe^t$, and assume that Theorem \ref{arHilbert}.5
%holds for subvarieties of codimension at most $s-1$.
%Then, $F_{D}(\CY)$ contains a non zero lattice point $f$ with
%\[ \log |g| \leq  \left(\frac{-c_1h(\CY_{s-1})}{\deg Y_{s-1}} + c_2 \right)
%   \frac{c}{(t-s+1)!} D.  \]
%\end{Satz}

\satz{Lemma} \label{mlvier}
Assume that Theorem \ref{arHilbert}.3
holds for subvarieties of codimension at most $s-1$, and let
$\CY \subset \Pe^t$ be an irreducible variety of codimension $s$.
With $\bar{D} = TR(\CY)$ the top regularity of $\CY$, and
$e_0= \frac1m$ there is a
$\tilde{D} \in \N$ with 
$D_0:=[e_0 \bar{D}] \leq \tilde{D} \leq 2\bar{D}$ such that
\[ \hat{H}_{\CY}(\tilde{D}) \geq (c_6 h(\CY) - c_6' \deg Y) \tilde{D}^{t+1-s}, \]
where $c_6,c_6'$ are constants only depending on $s$ and $t$.
\end{Satz}

\proof
If
\[ \hat{H}_{\CY}(\bar{D}) \geq e_2(h(\CY) - \deg Y) \bar{D}^{t+1-s}, \]
choose $\tilde{D} = \bar{D}, c_6=c_6'=e_2$. If not,
let $\bar{\CY},\CY_{s-1},D_0$ be as in Corollary \ref{Ysm1}.
With $C=8e_5$ we make the case distinction
$h(\CY) \leq C h(\bar{\CY}) \bar{D}^{s-r}$ or 
$h(\CY)>C h(\bar{\CY}) \bar{D}^{s-r}$.
If $h(\CY) \leq C h(\bar{\CY}) \bar{D}^{s-r}$, 
since $TR(\bar{Y}) \leq D_r \frac{D_0}{2t}\leq \frac{D_0+1}2$, 
and Theorem \ref{arHilbert}.3 is assumed to hold for
subvarietes of codimension at most $s-1$, Corollary \ref{Ysm1}.1 implies
\[ \hat{H}_{\CY}(D_0+1) = \hat{H}_{\bar{\CY}}(D_0+1) \geq
   (c_3(r,t) h(\bar{\CY}) - c_3'(r,t) \deg \bar{Y}) (D_0+1)^{t+1-r}, \]
which by assumption and Corollary \ref{Ysm1}.3  is greater or equal
\[ \frac{c_3(r,t)}C h(\CY) \bar{D}^{r-s} (D_0+1)^{t+1-r} - 
   \frac{c_3'(r,t)}{e_3}\deg Y \bar{D}^{r-s} (D_0+1)^{t+1-r} \geq \]
\[ \left(\frac{c_3(r,t)}C h(\CY) - 
   \frac{c_3'(r,t)e_0^{r-s}}{e_3}\deg Y\right) D_0^{t+1-s}. \]
Choosing $\tilde{D}=D_0+1,c_6(s,t)= 
\max_{i<s}\frac{c_3(i,t)}C,c_6(s,t)'= \max_{i<s}\frac{c_3'(i,t)e_0^{r-s}}{e_3}$, the 
claim follows.

If $h(\CY) > C h(\bar{\CY}) \bar{D}^{s-r}$, by Corollary \ref{Ysm1}.3,
\[ h(\CY_{s-1}) \leq \frac18 \frac{h(\CY)}{\bar{D}} + 
   \frac{e_5}C \frac{h(\CY)}{\bar{D}}  + \frac 18 \frac{\deg Y}{\bar{D}}  \leq 
   \frac14 \frac{h(\CY)}{\bar{D}} + \frac18 \frac{\deg Y}{\bar{D}}.  \]
With $\tilde{D} =\left[\left(1+\frac1m\right)\bar{D}\right] \leq 2 \bar{D}$, 
by the previous Lemma, for every  vector 
$g \in I_\CY(\tilde{D}) \setminus I_{\bar{\CY}_{s-1}}(\tilde{D})$,
\[ \log |g^\bot_{Y_{s-1}}|_{L^2(\Pe^t)} \geq \frac{h(\CY)}{\deg Y_{s-1}} - 
   \tilde{D} \frac{h(\CY_{s-1})}{\deg Y_{s-1}}- c_5\tilde{D} \geq \]
\[ \frac{h(\CY)}{\deg Y_{s-1}} 
   -\frac{\tilde{D}}{4\bar{D}}\frac{h(\CY)}{\deg Y_{s-1}}
   - \frac{\tilde{D}}{8\bar{D}}\frac {\deg Y}{\deg Y_{s-1}}- c_5\tilde{D}\geq 
   \frac{h(\CY)}{2 \deg Y_{s-1}}  - (c_5+1) \tilde{D}, \]
where we used Corollary \ref{Ysm1} to estimate $h(\CY_{s-1})$.
Further, by Corollary \ref{Ysm1}.1
\[ \rk I_Y(\tilde{D})/I_{Y_{s-1}}(\tilde{D}) \geq e_2 \deg Y \tilde{D}^{t-s} \geq 
   \frac{e_2 e_4}2 \deg Y_{s-1} \tilde{D}^{t+1-s}. \]
Hence by the Theorem of Minkowski and Lemma \ref{mlnull},
\begin{eqnarray*}
\hat{H}_Y (\tilde{D}) &\geq& - \wde \; \bar{I}_Y(\tilde{D}) \geq
   -\wde \bar{I}_Y(\tilde{D})/\bar{I}_{Y_{s-1}}(\tilde{D}) \\
   &\geq& \frac{e_2 e_4}2 \deg Y_{s-1} \tilde{D}^{t+1-s}
   \left(\frac{h(\CY)}{2 \deg Y_{s-1}}  - (c_5+1)\tilde{D} \right)-
   H_Y(\tilde{D}) \log H_Y(\tilde{D}) \\
  & \geq&  \frac{e_2 e_4}4 h(\CY) \tilde{D}^{t+1-s}- 
   \frac{(c_5+1) e_2 e_4}2 \deg Y_{s-1} \tilde{D}^{t+2-s} \\
  &-&  \deg Y (2t\tilde{D})^{t-s} \log \left(\deg Y (2t\tilde{D})^{t-s} \right)\\
  &\geq& \frac{e_2 e_4}4 h(\CY) \tilde{D}^{t+1-s}- 
   (c_5+1) e_2 \deg Y \tilde{D}^{t+1-s}-
   c' \deg Y \tilde{D}^{t-s} \log \tilde{D}, 
\end{eqnarray*} 
with a suitable $c'$ only depending on $t$ and $s$. The claim thus follows with
$c_6 = \frac{e_2e_4}4, c_6'=(c_5+1)e_2+c'$,

\proof {\sc of Theorem \ref{arHilbert}.4:}
The proof is by complete induction on the codimension of $\CY$. If
the codimension is $0$, then $\CY = \Pe^t$, and by Lemma \ref{basisnorm},
\[ \hat{H}_{\CX}(\tilde{D}) \geq \frac{\sigma_t}2 
   \frac{\tilde{D}^{t+1}}{(t+1)!}, \]
if $D \geq 0 = \bar{D}$.
Assume now the Theorem is proved for subvarietes of codimension at most
$s-1$, and let $\bar{D}$ again be the top regularity of $\CY$, and
$\tilde{D}$ as in the previous Lemma.
For $2 \bar{D} \leq D \leq 4t \bar{D}$, with $\tilde{D}$, as $\CY$ is 
irreducible Lemma, \ref{mlnull} and Proposition \ref{alHilbert}.1 imply
\[ \hat{H}_{\CY}(D) \geq \hat{H}_{\CY}(\tilde{D}) -c_1'D H_Y(\tilde{D})
\geq 
   (c_6 h(\CY)-c_6'\deg Y)\tilde{D}^{t+1-s}- c_1't^t\deg Y D \tilde{D}^{t-s} 
  \geq \]
\begin{equation} \label{fall1}
(c_6 h(\CY)-(c_6'+c_1't^t)\deg Y) \left(\frac1{(4t)^2}\right)^{t+1-s} D^{t+1-s}. 
\end{equation}

For $D \geq 4t \bar{D}$, since by the previous Lemma,
\[ \hat{H}_{\CY}(\tilde{D}) \geq (c_6 h(\CY)-c_6'\deg Y)\tilde{D}^{t+1-s}, \]
and by Theorem \ref{alHilbert}.1,
\[ H_Y(\tilde{D}) \leq \deg Y {\tilde{D}+t-s \choose t-s} \leq
   \deg Y (t \tilde{D})^{t-s} . \]
the Theorem of Minkowski assures the existence of a nonzero vector
$f \in F_{\CY}(\tilde{D})$ with
\begin{eqnarray*}
 \log |f| &\leq&  -\frac{c_6h(\CY)-c_6'\deg Y}{t^{t-s}\deg Y} \tilde{D} + 
   \frac12 \log H_Y(\tilde{D})   \\
  &\leq& -\frac{c_6 h(\CY)-c_6' \deg Y}{t^{t-s}\deg Y} \tilde{D} + 
   \frac12 \log \left(\deg Y t^{t-s} \tilde{D}^{t-s}\right) \\
  &\leq& -\frac{c_6}{t^{t-s}} \frac{h(\CY)}{\deg Y} \tilde{D} + c' \tilde{D}, 
\end{eqnarray*}
with a positive $c'$ suitably chosen.

Put 
$l = \left[\frac{D}{2\tilde{D}}\right]\geq \frac{4\bar{D}}{2\tilde{D}}\geq 1$ 
and
$\bar{f} := f^{\otimes l} \in \Gamma(\CY,O(l\tilde{D}))$. By Lemma \ref{polprod},
\[ \log |\bar{f}| \leq l \log |f| + c_1 l \tilde{D} \leq
   - \frac{c_6h(\CY)}{t^{t-s}\deg Y} l\tilde{D} + (c'+ c_1) l\tilde{D}. \]
Thus, by Lemma \ref{mlnull}, 
the lattice $\bar{f}F_Y(D-l\tilde{D}) \subset F_Y(D)$ 
is spanned by vectors of logarithmic length at most 
$\log |\bar{f}| \leq  -\frac{c_6}{t^{t-s}} 
\frac{h(\CY)}{\deg Y} l \tilde{D} + (c'+ c_1) l\tilde{D}$. 
Further,
\[ \rk \; \left(\bar{f} F_Y(D-l\tilde{D})\right) = 
   \rk \; F_Y(D-l\tilde{D})\geq H_Y(D-l\tilde{D}). \]
As $l\tilde{D} \leq \frac{D}2$, one deduces
$D -l\tilde{D} \geq \left[\frac D2\right]$, hence the above is greater or equal
$H_Y\left(\left[\frac{D}2\right]\right)$ which using Theorem 
\ref{alHilbert} and the inequality 
$\left[\frac D2\right]-s\bar{D} \geq \left[\frac D4\right]$
is greater or equal
\[ \deg Y {\left[\frac{D}2\right] - sD +t-s  \choose t-s} \geq
   \deg Y \frac{\left[\frac{D}4\right]^{t-s}}{(t-s)!} \geq
   \frac{\deg Y}{8^{t-s}(t-s)!} D^{t-s}. \]
Thus Lemma \ref{mlnull} implies
\begin{eqnarray*}
 \hat{H}_Y(D)  &\geq& 
   \wde \overline{\left(\bar{f} F_Y(D-l\tilde{D})\right)} 
   \geq - \rk \; \left(\bar{f} F_Y(D-l\tilde{D})\right) \log |\bar{f}| \\
 &\geq&  \frac{c_6}{(8t)^{t-s}(t-s)!} h(\CY) l \tilde{D} D^{t-s} -
   \frac{c'+c_1}{8^{t-s}(t-s)!} \deg Y l\tilde{D} D^{t-s}. 
\end{eqnarray*}
As $\tilde{D} \leq (1 + \frac1{2t}) \bar{D}$, and 
$\bar{D} \leq \frac D{4t}$, one concludes $l \tilde{D} \geq \frac D4$. 
Since also $l \tilde{D} \leq D$, one arrives at
\begin{equation} \label{fall2}
   \hat{H}_Y(D) \geq \frac{c_6}{4^{t+1-s}(t-s)!} 
   h(\CY)D^{t+1-s} -
   \frac{c'+c_1}{4^{t-s}(t-s)!} \deg Y D^{t+1-s}. 
\end{equation}
In view of (\ref{fall1}) and (\ref{fall2}), setting
\[ c_3(s) = \max \left(\frac{c_6(s-1)}{(4t^2)^{t+1-s}},  
   \frac{c_6(s-1)}{8^{t+1-s}(t-s)!}\right), \]
\[ c_3'(s) = \max \left(\frac{c_6'(s-1)}{((4t)^2)^{t+1-s}}+c_1'(2t)^t,
   \frac{c'+c_1}{8^{t-s}(t-s)!}\right). \]
finishes the proof.

\begin{appendix}

\section{A simpler proof of a weakened version of Theorem \ref{alHilbert}}

\satz{Proposition}
Let $k$ be a field of characteristic $0$, 
and $X$ an irreducible subvariety of codimension $s$ of $\Pe^t_k$.
There is a constant $c_1$ only depending on $s$ and $t$ such
that for $D \geq (t+1-s)\bar{D} = (t-s) TR(X)+1$ the inequality
\[ H_X(D) \geq c_1(s,t)D^{t-s} \]
holds.
\end{Satz}

\vspace{2mm}

Before going into the proof, I collect a few well known facts on projective
and affine schemes over a field.
There are the canonical $k$-linear maps
\[ \varphi_0: k[x_0,\ldots,x_t] \to k[x_1,\ldots,x_t], \quad
             f(x_0,x_1,\ldots,x_t) \mapsto f(1,x_1,\ldots,x_t), \]
\[ \psi_1:  k[x_1,\ldots,x_t] \to k(x_1)[x_2,\ldots,x_t], \quad
            f(x_1,\ldots,x_t) \mapsto f(x_1,\ldots,x_t). \]

\satz{Lemma} \label{aeins}
\begin{enumerate}
\item
With $U_0 = \Pe^t \setminus V(x_0)$ the induced map 
$(\varphi_0)^*: \A^t \to U_0$, is a bijection from the set of 
prime ideals of $k[x_1,\ldots,x_t]$ to 
the set of homogeneous prime ideals of $k[x_0,\ldots,x_t]$ that do not
contain $x_0$ with inverse $(\varphi_0^*)^{-1}: \fP \mapsto \varphi_0(\fP)$. 
Moreover, $(\varphi_0^*)^{-1}$ commutes with set theoretic operations.
\item
If $\A^t_{x_1}$ is defined as the set of prime ideals $\fP$ in 
$k[x_1,\ldots,x_n]$ such that
$\fP \cap k[x_1] = \{0\}$, then $x \in \A^t_{x_1}$ if and only if
the map $X= \bar{x} \hookrightarrow \A^t \to \A^1$ is dominant and
$\psi_1^*:\A^{t-1}_{k(x_1)} \to \A^t_{x_1}$ is a bijection with inverse
$(\psi_1^*)^{-1}: \fP \mapsto (\psi_1(\fP))$.
$\A^{t-1}_{k(x_1)}$. Moreover, $(\psi_1^*)^{-1}$ commutes with set theoretic 
operations. 
\end{enumerate}
\end{Satz}

\satz{Lemma} \label{azwei}
With the notations of the previous Lemma, let $x \in U_0 \subset \Pe^t$ 
be a homogeneous prime ideal of dimension $p$ such that 
for every homogeneous polynomial $h \in k[x_0,x_1]$, the intersection of
$X=\bar{x}$ with $V(h) \cap V(x_2) \cap \cdots \cap V(x_p)$ is proper,
hence $x^\circ=(\varphi_0^*)^{-1}(x) \in \A_{x_1}^t$. Put
$y^\circ = (\psi_1^*)^{-1} (x^\circ)$, and $y= (\varphi')_0^*(y^\circ)$,
where $\varphi'_0$ is the map from $k(x_1)[x_0,x_2,\ldots,x_t]$ to
$k(x_1)[x_2,\ldots,x_t]$ that maps $x_0$ to $1$, and let
$X= \bar{x},X^\circ = \bar{x^\circ},Y^\circ = \bar{y^\circ},Y=\bar{y}$. Then,
\begin{enumerate}
\item
\[ \mbox{TR} (Y) \leq \mbox{TR}(X). \]
\item
\[ \deg Y = \deg X. \]
\item
For every $D \in \N$,
\[ H_X((p+1)D) \geq D H_Y(pD). \]
\end{enumerate}
\end{Satz}

\proof
1. Let $X = \overline{U \cap V(f_1) \cap \cdots \cap V(f_{t-p})}$ with
$\deg f_i \leq \mbox{TR}(X), i=1,\ldots,t-p$. 
Since $X$ intersects $V(x_0)$ properly, one may replace $U$ by 
$U'= U \cap \A^t$. By Lemma \ref{aeins}.1, 
\[ X^\circ = \overline{(\varphi_0^*)^{-1} x}  =
   \overline{U' \cap V(\varphi_0(f_1)) \cap \cdots V(\varphi_0(f_{t-p}))}. \]
Further, since $x$ does not contain $x_0$, it may be assumed
that each $f_i$ is not a multiple of $x_0$, hence 
$\deg f_i = \deg \varphi_0(f_i)$.

Next, by Lemma \ref{aeins}.2, with $U''= (\psi_1^*)^{-1}(U' \cap \A^t_{x_1})$,
\[ Y^\circ = \overline{(\psi_1^*)^{-1} x^{\circ}} = 
   \overline{U'' \cap V((\psi_1\circ \varphi_0) f_1) \cap \cdots \cap
               V((\psi_1\circ \varphi_0) f_{t-p})}, \]
and consequently, 
\[ Y = \overline{\varphi_o^{-1}(U'') 
   \cap V((\varphi_0^{-1}\circ\psi_1\circ \varphi_0) f_1) \cap 
                    \cdots \cap
               V((\varphi_0^{-1}\circ\psi_1\circ \varphi_0) f_{t-p})}. \]
Since $\deg_{k(x_1)} f$ is the degree of $f$ in the variables
$x_2,\ldots,x_t$, one gets
\[ \deg_{k(x_1)} \varphi_0^{-1}\circ \psi_1\circ \varphi_0 f_i =
   \deg_{k(x_1)} \psi_1\circ \varphi_0 f_i  \leq
   \deg_k \varphi_0 f_i = \deg _k f_i \leq \mbox{TR}(X), \]
for every $i=1,\ldots,t-p$, finishing the proof.

\vspace{2mm}

2. Put $\varphi_0':k(x_1)[x_0,x_2\ldots,x_t]\to k(x_1)[x_2,\ldots,x_t], \;\;
f(x_0,\ldots,x_t) \mapsto f(1,x_1,\ldots,x_t)$.

Since $X \cap V(x_0) \cap \cdots \cap V(x_p) = \emptyset$, 
for every $i=p+1,\ldots,t$ there is an $l_i\in \N$ such that
$x_i^{l_i} \in x + (x_0,\ldots,x_p)$, hence there is an $f_i \in x$ as well
as $g_{ij} \in k[x_0,\ldots,x_t], j=0,\ldots,p$, homogeneous polynomials
such that
\[ x_i^{l_i} = \sum_{j=0}^p x_j g_{ij} + f_i. \]
Hence,
\[ x_i^{l_i} = x_0 \bar{g}_{i0} + x_0 x_1 \bar{g}_{i1}+
   \sum_{j=2}^p x_j \bar{g}_{ij} + \bar{f}_i, \]
where $\bar{g}_{ij},\; i=p+1,\ldots,t,j=1,\ldots,p$ are homogeneous polynomials 
that $\varphi_0'(\bar{g}_{ij}) = \psi_1 \circ \varphi_0(g_{ij})$, and
$\bar{f}_i, \; i=p+1,\ldots,t$ are homogeneous polynomials with
$\varphi_0'(\bar{f}_i) = \psi_1 \circ \varphi_0(f_i)$.

%Since $((\varphi_0')^{-1} \circ \psi_1 \circ \varphi_0) (x_0^{j_0} \cdots x_t^{%j_t})=
%x_1^{j_1} x_0^{j_0+j_1} x_2^{j_2} \cdots x_t^{j_t}$ for every monomial, this
%implies

%\[ (x_i^{l_i}) = x_0 ((\varphi_0')^{-1} \circ \psi_1 \circ \varphi_0) (g_{i0}) +
%      x_1 x_0 ((\varphi_0')^{-1} \circ \psi_1 \circ \varphi_0) (g_{i1}) + \]
%\[ \sum_{j=2}^p x_j ((\varphi_0')^{-1} \circ \psi_1 \circ \varphi_0) (g_{ij})´+
%     ((\varphi_0')^{-1} \circ \psi_1 \circ \varphi_0) (f_i), \]
%for every $i=p+1,\ldots,t$, 
Thus,
$x_i^{l_i} \in  (((\varphi_0')^*)^{-1} 
\circ (\psi_1^*)^{-1} \circ (\varphi_0^*)^{-1})x+
(x_0,x_2,\ldots,x_p)$, i.e.\@

$Y \cap V(x_0) \cap V(x_2) \cap \cdots \cap V(x_p) = \emptyset$, and thus
\[ \deg Y = \deg Y_1 \]
with $Y_1:=Y . \di x_2 . \ldots . x_t=Y^\circ . \di x_2 . \ldots . x_t$,
It remains to prove $\deg X = \deg Y_1$.

Since for any linear form $L(x_0,x_1)= a x_0+bx_1, \; a,b \in k$, 
the intersection of \\
$X_1=X . \di x_2 . \ldots . \di x_p$ with $\di L$ is proper, the degree
of $X_1 . \di L$ equals the degree of $X$ and also the degree of 
the map $\varphi$, which in turn equals the degree of the curve
$X_1^\circ$ over $k[x_1]$ which by definition is the degree of $Y_1$

\vspace{2mm}

3. Firstly, $\varphi_0,\varphi_0'$ map
$\Gamma(\Pe^t_k,O(D)) \cong k[x_0,\ldots,x_t]_{D,hom}$ and 
$\Gamma(\Pe^{t-1}_{k(x)},O(D)) \cong k(x_1)[x_0,x_2,\ldots,x_t]_{D,hom}$ 
isomorphically to $k[x_1,\ldots,x_t]_{\leq D}$ and
$k(x_1)[x_2,\ldots,x_t]_{\leq D}$ respectively. Since they also identify
$I_X(D)$ with $I_{X^\circ}(D)$ and $I_Y(D)$ with
$I_{Y^\circ}(D)$, it suffices to prove
\[ \rk \; \frac{k[x_1,\ldots,x_t]_{\leq (p+1)D}}{I_{X^\circ}((p+1)D)} \geq
   D \; \rk \; \frac{k(x_1)[x_2,\ldots,x_t]_{\leq pD}}{I_{Y^\circ}(pD)}. \]

Let $A$ be the set of multi indices $I=(i_2,\ldots,i_t)$ with
$|I| := i_2 + \cdots i_t = pD$. Since $(x_2^{i_2} \cdots x_t^{i_t})_{I \in A}$ 
generate $k(x_1)[x_2,\ldots,x_t]_{\leq pD}$, there
is a subset $B \subset A$ such that $(x_2^{i_2} \cdots x_t^{i_t})_{I \in B}$
is a basis of $\frac{k(x_1)[x_2,\ldots,x_t]_{\leq pD}}{I_{Y^\circ}(pD)}$ over
$k(x_1)$.
I claim that
\[ \CB = \left\{ x_1^{i_1} x_2^{i_2} \cdots x_t^{i_t} |
           i_1 \leq D, (i_2,\ldots,i_t) \in B \right\} \]
are linearily independent over $k$ in 
$\frac{k[x_1,\ldots,x_t]_{\leq (p+1)D}}{I_{X^\circ}((p+1)D)}$ which, of course,
implies the above inequality. 

So let $a_{i_1,\ldots,i_t} \in k$ such that
\[ f=\sum_{i_1 \leq D, (i_2,\ldots,i_t) \in B} a_{i_1,\ldots,i_t} x_1^{i_1} \cdots x_t^{i_t} \in
   I_{X^\circ}((p+1)D). \]
Then,
\[ \psi_1(f) =
   \sum_{i_1 \leq D, (i_2,\ldots,i_t) \in B} a_{i_1,\ldots,i_t} x_1^{i_1} \cdots x_t^{i_t} \in
   I_{Y^\circ}(pD). \]
Since the $(x_2^{i_2} \cdots x_t^{i_t})_{I \in B}$ are linearily independent
over $k(x_1)$ in $\frac{k(x_1)[x_2,\ldots,x_t]_{\leq pD}}{I_{Y^\circ}(pD)}$, 
this implies $\sum_{i_1 \leq D} a_{i_1,\ldots,i_t} x_1^{i_1} \in I_{Y^\circ}(D)$, for
all $(i_2,\ldots,i_t)\in B$, hence
\[ \sum_{i_1 \leq D} a_{i_1,\ldots,i_t} x_1^{i_1} \in I_{Y^\circ}(D) \cap k[x_1] =
   \{0\}  \quad \forall (i_2,\ldots,i_t)\in B \]
by assumption, which implies
$a_{i_1,\ldots,i_t} = 0$ for all $i_1 \leq D, (i_2,\ldots,i_t) \in B$.

\proof {\sc of the Proposition}
Set $p=t-s$ the dimension of $X$. The proof will be given by complete 
induction over $t+s$. For $t+s=0$, i.e.\@ $s=t=0$ the Proposition trivially
holds. Assume the Proposition holds for $k$, let
$X$ be of codimension $s$ in $\Pe^t$ with $s+t=k+1$, and make the
case distinction $s=t$ or $s<t$.

If $s<t$ and thus $p\geq 1$, making a linear coordinate transformation if 
necessary, one may assume that $X$ does not intersect
$V(x_0)\cap \cdots \cap V(x_p)$, i.e.\@ with $x$ the homogeneous
prime ideal in $k[x_0,\ldots,x_t]$ such that $X=\bar{x}$, we have
$x \cap (x_0,\ldots,x_p) = \{0\}$, and
for every linear form $L \in \la x_0,x_1 \ra$ the intersection of
$X$ with $V(L) \cap V(x_2) \cap \cdots \cap V(x_p)$ is proper. 
With $X^\circ,Y^\circ,Y$ as in Lemma \ref{azwei},
$\mbox{TR}(Y) \leq \mbox{TR}(X)=\bar{D}$.

For $D \geq (p+1) \bar{D}$, put $D_1:=\left[\frac D{p+1}\right]$. 
Then, $pD_1 \geq p\bar{D}$, and by Lemma \ref{azwei}.3 and induction hypothesis,
\[ H_X(D) \geq H_X((p+1)D_1) \geq D_1 H_Y(pD_1) \geq \]
\[ D_1 c_1(t-1,s) \deg Y D_1^p = c_1(t-1,s) \deg X D_1^{p+1}, \]
the equality holding because of Lemma \ref{azwei}.2. 
Since $p \geq 1$, we have $D_1 \geq \frac D2$, hence the above is greater or
equal than
\[ \frac{c_1(t-1),s}{2^{p+1}} \deg X D^{p+1}, \]
proving the claim with $c_1(t,s) =  \frac{c_1(t-1),s}{2^{p+1}}$.

If $s=t$, set $m=2t$, $D_r = \left[\frac {\bar{D}}{m^{s+1-r}}\right]$,
and let $r$ be the maximal number such that $X$ is contained in $\bar{X}$ a 
locally complete intersection of $r$ hypersurfaces of degree at most $D_r$
at $X$. By Lemma \ref{abschn}, there is an irreducible component $\bar{Y}$ of
$\bar{X}$ such that $I_X(D_{r+1})\subset I_{\bar{Y}}(D_{r+1})$, and
$\deg X \leq \deg \bar{Y} \bar{D}^{s-r}= \deg \bar{Y} \bar{D}^{t-r}$. 
Since, $D_{r+1}+1 \geq mD_r \geq (t+1-r) D_r \geq (t+1-r) \mbox{TR}(\bar{Y})$
by induction hypothesis,
\[ H_X(D_{r+1}+1) = H_{\bar{Y}}(D_{r+1}+1) \geq 
   c_1(t,r) \deg \bar{Y} (D_{r+1}+1)^{t-r} \geq \]
\[ c_1(t,r) \deg X (D_{r+1}+1)^{t-r}\bar{D}^{r-t}. \]
Thus, for $D \geq (t+1-s) \bar{D}+1 = \bar{D}+1 \geq D_{r+1}+1$,
\[ H_X(D) \geq H_Y(D_{r+1}+1) \geq 
   \frac1{m^{(s+1-r)(t-r)}} \deg X D_{r+1}^{t-t} = \]
\[  \frac1{m^{(s+1-r)(t-r)}} \deg X D^{t-t}, \]
proving the claim with $c_1(t,t)= \frac1{m^{(s+1-r)(t-r)}}$.

\end{appendix}

\end{document}